\documentclass[11pt,a4paper,twoside,reqno]{amsart}
\usepackage{amsfonts,amssymb,amscd,amsmath,enumerate,verbatim,calc}

\usepackage{mathrsfs}

\usepackage{rotating}
\newcommand{\ud}{\, \textrm{d}}

\newcommand{\ui}{\, \indent\indent}
\makeindex
\usepackage{xcolor}

\begin{document}
\begin{center}{\Large{An expository note on Prohorov metric and Prohorov Theorem\\ by\\ R.P. Pakshirajan$^1$ and M. Sreehari$^2$ }}\\
	$^1$ 227, 18th Main, 6th Block, Koramangala\\
	Bengaluru- 560095, Karnataka, India.\\
 E-mail: \textit{vainatheyarajan@yahoo.in} \\
	$^2$ 6-B, Vrundavan Park, New Sama Road\\
	Vadodara- 390024, Gujarat, India.\\
	Corresponding author E-Mail: \textit{msreehari03@yahoo.co.uk}\\
\end{center}

\clearpage
\textbf{Abstract}
The main aim of this article  is to give an exposition of weak convergence, Prohorov theorem and Prohorov  spaces. In this context we study the relationship between  Levy distance $\ell(F, G)$ between two distribution functions $F$ and $G$  and the Prohorov distance $\pi( \mu, \nu)$  between the probability measures $\mu$ and $\nu$  determined by $F$ and $G$ respectively. We study the relationship among  the weak convergence of probability measures  ($\mu_n$)  determined by  distribution functions ($F_n$) to the probability measure  $\mu$ determined by a distribution function $G$, the convergence of  $\ell(F_n, G)$ and $\pi( \mu_n, \nu)$ to zero  under suitable assumptions on the metric space on which these measures are defined. Tightness of probability measures and  relative sequential compactness are studied and Prohorov theorem is proved in different settings. Prohorov  spaces and non-Prohorov spaces are discussed.
\clearpage

	\section{Introduction and Prohorov distance}
The aim of this expository article is to discuss some of the contributions of the  mathematician Yuri V. Prohorov to the field of probability theory. 	In particular we discuss (a) Prohorov distance function and its relation to Levy distance function on the space of probability measures on a separable metric space, (b) Prohorov theorem which has useful impact on the study of functional limit theorems/ stochastic processes, and (c) Prohorov space and how the absence of sufficient structures in the topological spaces where the measures are defined renders those spaces ineligible to be a Prohorov space.\\
Before we define the two distance functions proposed by Prohorov and Levy, we may mention that in probabilty theory various distance functions are proposed and their properties are studied. We refer to Gibbs and Su [3] for details. \\
Before we define the Prohorov distance function we need to introduce $\varepsilon$-neighborhood of a set and prove a related result needed later. \\\\
\noindent \textbf{Definition 1.1.} \\
In a metric space $(\textbf{M},\ d)$, by the $\varepsilon$-neighborhood  $A^{\varepsilon}$\ of a subset $A$, we understand the set $\{x \in \textbf{M}\ \big| \;  \exists\; y \in A\ \text{with}\ d(x,\ y) < \varepsilon \}.$
It is easy to see \\\ui 
(i) that $A^{\varepsilon} = \bar{A}^{^{^{\varepsilon}}}$, $\bar{A}$\ being the closure of $A$,\ and \\\ui 
(ii) that, if $A$\ is a closed set then  $A^{\frac{1}{r}} \downarrow A$ as $r\uparrow \infty$.\\
\noindent \textbf{Theorem  1.1.}\\\hspace*{.5cm}For a sequence $(\mu_n),\ n = 0,\ 1,\ 2,\ \ldots$\ of probability measures on  $\textit{\Large{m}}$, the Borel $\sigma$-field of  $\textbf{M}$,  the following four conditions are equivalent:\\
(i)\hspace*{2.6cm}$\varlimsup\limits_{n \rightarrow \infty}\mu_n(C) \le \mu_0(C)$\ for every closed set $C$.\\(ii)\hspace*{2.5cm}$\varliminf\limits_{n \rightarrow \infty}\mu_n(D) \ge \mu_0(D)$\ for all open sets $D$.\\(iii)\hspace*{2.5cm}$\lim\limits_{n \rightarrow \infty}\mu_n(A) = \mu_0(A)$\ for every set $A$\ with $\mu_0(\partial A) = 0$\\\hspace*{5.5cm}where $\partial A$\ stands for the boundary of $A$.\\
(iv)\hspace*{2.5cm}$\lim\limits_{n \rightarrow \infty}\int\limits_{\textbf{M}} f(x)\ud\mu_n(x) = \int\limits_{\textbf{M}}f(x)\ud\mu(x)$\ for every real bounded\\\hspace*{7cm} uniformly continuous function $f$\\\hspace*{7cm}satisfying $0 \le f(x) \le 1,\ x \in \textbf{M}$.\\
This is a wellknown result and we refer to, for example, Theorem 2.1, p. 16, Billingsley [1]; Theorem 2.4.2, p. 98 in Pakshirajan [7]\\
We now introduce two important concepts in probability and relate them to the main interest of this Section, namely the distance functions introdued by Levy and Prohorov.\\
\textbf{Definition 1.2.}\\
A sequence $(\mu_n),\ n = 0,\ 1,\ 2,\ ...$\ of probability measures on $\textit{\Large{m}}$\ is said to converge weakly to $\mu_0\ (\mu_n\ \ ^{\underrightarrow{w}}\ \ \mu_0)$\ \\\hspace*{2.5cm} if (i) or (ii) or (iii) or (iv)   of Theorem 1.1 \; holds.\\
\textbf{Definition 1.3.}\\ A function $F: R\rightarrow[0, 1]$ is called a distribution function if $F$ is monotone, right continuous, $F(-\infty) =0$ and $F(+\infty)=1.$\\
Let $(\textbf{M},\ d)$\ be a metric space and let $\mathcal{M}$\ denote the totality of all the probability measures on {\textit{\Large{m}}}. When $\textbf{M}$\ is the real line $\textbf{R}$, the totality of all the distribution functions on $\textbf{R}$\ will be denoted by $\mathcal{F^*}$\\
\textbf{Definition 1.4.}\\
For $F,\ G \in \mathcal{F^*}$,  the Levy  distance $\ell$ is given by \\$\ell(F,\ G) = \inf\{h > 0 \big| F(x - h) - h\le G(x)\le F(x + h) + h\ \text{for all}\ x \in \textbf{R}\}.$(1.1)	
\textbf{Definition 1.5.}\\
For $\mu,\ \nu \in \mathcal{M}$, the Prohorov distance $\pi$ is given by \\ $\pi(\mu,\ \nu) = \inf\{\varepsilon > 0 \big| \mu(A) \le \nu(A^{\varepsilon}) + \varepsilon\ \text{and}\ \nu(A) \le \mu(A^{\varepsilon}) + \varepsilon\\\hspace*{7.5cm} \text{for all closed sets}\ A \}.$\hfill(1.2)
Trivially, $0 \le \ell(F,\ G) \le 1,\ F,\ G \in \mathcal{F^*};\ 0 \le \pi(\mu,\ \nu) \le 1,\ \mu,\ \nu \in \mathcal{M}.$\\
Let us denote the probability measures determined by the two distribution functions $F$ and $G$ by $\mu$ and $\nu$ respectively.\\
For completeness we consider a related distance function and prove a related result useful later.\\
$\pi^*(\mu,\ \nu) = \inf\{\varepsilon > 0 \big| \mu(A) \le \nu(A^{\varepsilon}) + \varepsilon\ \text{and}\ \nu(A) \le \mu(A^{\varepsilon}) + \varepsilon\\\hspace*{7.5cm} \text{for all sets}\ A \in \textit {\Large{m}} \}.$\hfill(1.3)
\noindent\textbf{Theorem 1.2.}\\
(i) $\pi^*(\mu,\ \nu) = \pi(\mu,\ \nu),\ \mu,\ \nu \in \mathcal{M}.$\\
(ii) If $\mu(E) \le \nu(E^{\varepsilon}) + \varepsilon\ \text{for every closed set}\ E$, then\\\hspace*{3cm} $\nu(E) \le \mu(E^{\varepsilon}) + \varepsilon\ \text{for every closed set}\ E$ and conversely.\\
Proof.
(i) Write $\pi(\mu,\ \nu) = \inf\{\varepsilon > 0 \big| \varepsilon \in \Lambda_1\}$\ and $\pi^*(\mu,\ \nu) = \inf\{\varepsilon > 0 \big| \varepsilon \in \Lambda_2\}$\ where $\Lambda_1=\{\varepsilon > 0 \big| \mu(A) \le \nu(A^{\varepsilon}) + \varepsilon\ \text{and}\ \nu(A) \le \mu(A^{\varepsilon}) + \varepsilon\; \text{for all closed sets}\ A \}$ and  $\Lambda_2=\{\varepsilon > 0 \big| \mu(A) \le \nu(A^{\varepsilon}) + \varepsilon\ \text{and}\ \nu(A) \le \mu(A^{\varepsilon}) + \varepsilon\; \text{for all sets}\ A \in \textit {\Large{m}} \} $. Note that $\Lambda_2 \subset \Lambda_1$. Hence $\pi(\mu,\ \nu) \le \pi^*(\mu,\ \nu)$. Now let $A \in \textit{\Large{m}}$\ be arbitrary and let $\varepsilon \in \Lambda_1$. Hence $\mu(A) \le \mu(\bar{A}) \le \nu({\bar{A}}^{^{^{\varepsilon}}}) + \varepsilon \le \nu(A^{\varepsilon}) + \varepsilon$\ and similarly,\ $\nu(A) \le \nu(\bar{A}) \le \mu({\bar{A}}^{^{^{\varepsilon}}}) + \varepsilon \le \mu(A^{\varepsilon}) + \varepsilon$ , thus showing that $\varepsilon \in \Lambda_2$. It now follows that $\Lambda_1 = \Lambda_2$. Hence $\pi(\mu,\ \nu) = \pi^*(\mu,\ \nu)$.\\
(ii) Let $\mu(E) \le \nu(E^{\varepsilon}) + \varepsilon$\ for every closed set {$E$}. We note $\Big(((E^{\varepsilon})^{\prime}\big)^{\varepsilon}\Big)^{\prime}$\ is a closed set. We further note, from the definition of $\varepsilon$-neighborhood, that $(E^{\varepsilon})^{\prime} \cap E = \emptyset$. It is also true that\\\hspace*{4cm} $\big((E^{\varepsilon})^{\prime}\big)^ {\varepsilon} \cap E = \emptyset.$\hfill(1.4)\\ For,  let $x \in \big((E^{\varepsilon})^{\prime}\big)^{\varepsilon}$. That would imply either (i) $x \in (E^{\varepsilon})^{\prime}$\ or (ii)  $x \notin (E^{\varepsilon})^{\prime} \ \text{and}\ d(x,\ y) \ge \varepsilon$\ for all $y \in (E^{\varepsilon}) ^{\prime}$. If (i) holds then 
$x {\notin}\; E^{\varepsilon}$. Hence $ x\notin\; E$. If (ii) holds, then there is a contradiction because it is possible that $x\;\in E^\varepsilon\sim E$ and $d(x, y) <\epsilon$ for some $y \in (E^\epsilon)^\prime$. Thus (1.4) is true and we get $E \subset \Big(((E^{\varepsilon})^{\prime}\big)^ {\varepsilon}\Big)^{\prime}$.\\ Hence
$\nu(E) \le \nu\Big(((E^{\varepsilon})^{\prime}\big)^ {\varepsilon}\Big)^{\prime}.$
Since $(E^{\varepsilon})^{\prime}$\ is a closed set, we have, by the hypothesis,\\\hspace*{3cm}
$\mu\big((E^{\varepsilon})^{\prime}\big) \le \nu\Big(((E^{\varepsilon})^{\prime}\big)^ {\varepsilon}\Big) + \varepsilon.$\\
Hence \hspace*{1cm}$\nu(E) \le \nu\Bigg(\Big(((E^{\varepsilon})^{\prime}\big)^ {\varepsilon}\Big)^{\prime}\Bigg) = 1 - \nu\Big(((E^{\varepsilon})^{\prime}\big)^ {\varepsilon}\Big) \le \mu(E^{\varepsilon}) + \varepsilon.$\qed\\
\textbf{Remark}.
By (ii) it follows that\\ $\pi(\mu,\ \nu) = \inf\{\varepsilon > 0 \big|\mu(A) \le \nu(A^{\varepsilon}) + \varepsilon,\; \text{for all closed sets}\ A \}. $\hfill{(1.5)}\\
Similarly\\ $\ell(F,\ G) = \inf\{\varepsilon > 0 \big| F(x) \le G(x + \varepsilon) + \varepsilon,\  \text{for all}\ x \in R\}.$ \hfill(1.6)\\
The following result establishes an important property of the two distance functions.\\
\noindent\textbf{Theorem 1.3.}\\
\hspace*{.5cm}  $\ell,\ \pi$\ \text{are proper metric functions}. \\
Proof.
We will prove the assertion only for $\pi$\ since the proof for $\ell$\ can be constructed on similar lines. 
We must show (i) $\pi(\mu,\ \nu) = \pi(\nu,\ \mu)$, (ii) $\pi(\mu,\ \nu) = 0$\ if and only if $\mu \equiv \nu$\ and (iii) the triangle inequality.
(i) By the very definition in (1.2), this symmetry property is assured.\\
(ii) If $\mu \equiv \nu$, then $\mu(A) = \nu(A)$\ for all closed sets $A$. Hence the inequalities  $\mu(A) \le \nu(A^{\varepsilon}) + \varepsilon\ \text{and}\ \nu(A) \le \mu(A^{\varepsilon}) + \varepsilon$\ hold for all closed sets $A$\ and all $\varepsilon$. Hence, by the definition of $\pi,\ \pi(\mu,\ \nu) = 0.$
Conversely, suppose $\pi(\mu,\ \nu) = 0$.Then for closed sets $A$, $\mu(A) \le \nu(A^{\varepsilon}) + \varepsilon$.  Let $\varepsilon \downarrow 0$ to get $\mu(A) \leq \nu(A)$. Similarly we have $\nu(A) \leq \mu(A)$.   The equality $\mu(A) = \nu(A)$\ for all closed sets implies, by the inner regularity property of measures in metric spaces (ref. p. 42, Pakshirajan [7]), that $\mu = \nu$
\ on \textit{\Large{m}}. 
(iii) We establish now the triangle inequality. First let us note the following. Suppose $\mu,\ \nu \in \mathcal{M}$. Suppose for $\varepsilon > 0$\ fixed, $\mu(A) \le \nu(A^{\varepsilon})$\ for every closed set $A$. Then for $\eta > 0$\\\hspace*{1cm} $\mu(A^{\eta}) = \sup\limits_{\substack{B \subset A^{\eta}\\ B\ \text{closed}}}\mu(B) \le \sup\limits_{\substack{B \subset A^{\eta}\\ B\ \text{closed}}} \nu(B^{\varepsilon}) \le \sup\limits_{\substack{B \subset A^{\eta}\\ B\ \text{closed}}}\; \sup\limits_{\substack{C \subset B^{\varepsilon}\\ C\ \text{closed} }}\nu(C)\\\\\hspace*{4.6cm} \le \sup\limits_{\substack{C \subset A^{\varepsilon + \eta}\\ C\ \text{closed} }}\nu(C) \le \nu(A^{\varepsilon + \eta}).$\\ 
Let $\mu_i \in \mathcal{M},\ i = 1,\ 2,\ 3$. Fix target error $\eta > 0$. Let\\{$\alpha \geq \pi(\mu_1,\;\mu_2)= \inf\{\varepsilon \big| \mu_1(A) \le \mu_2(A^{\varepsilon}) + \varepsilon\  \text{for all closed sets}\ A\}$ \\ {such} that $0 \le \alpha - \pi(\mu_1,\ \mu_2) < \eta$. Let \\
	$\beta \geq \pi(\mu_2,\;\mu_3)= \inf\{\varepsilon \big| \mu_2(A) \le \mu_3(A^{\varepsilon}) + \varepsilon\ \text{for all closed sets}\ A\}$\\ {such} that $0 \le \beta - \pi(\mu_2,\ \mu_3) < \eta$. 
	Let $A$\ be an arbitrary closed set. From $\mu_1(A) \le \mu_2(A^{\alpha}) + \alpha,\ \mu_2(A) \le \mu_1(A^{\alpha}) + \alpha$; $\mu_2(A) \le \mu_3(A^{\beta}) + \beta\ \text{and}\ \mu_3(A) \le \mu_2(A^{\beta}) + \beta$, we get :
	$\mu_1(A) \le \mu_2(A^{\alpha}) + \alpha \le \mu_3(A^{\alpha + \beta}) + \alpha + \beta$;\\
	\hspace*{3.7cm} again, $\mu_3(A)  \le \mu_2(A^{\beta}) + \beta \le \mu_1(A^{\alpha + \beta}) + \alpha + \beta$. \\ Hence $\pi(\mu_1,\ \mu_3) \le  \alpha + \beta \le \pi(\mu_1,\ \mu_2) + \pi(\mu_2,\ \mu_3) + 2\eta$. Since $\eta > 0$\ is arbitrary, the proof of the triangle inequality is complete.\qed\\\\
	We shall now discuss an interesting example to demonstrate how the Levy and Prohorov metrics can be computed in specific situations. The example also throws light on possible relations between them.\\\\
		\noindent\textbf{An example}.
	Let us calculate $\pi(\mu,\ \nu)$\ where $\mu,\ \nu$\ are measures generated by the distribution functions $F,\ G$:
	\begin{equation*}
	F(x) =
	\begin{cases}
	0 & \text{if}\ x < 0\\
	x & \text{if}\ 0 \le x \le 1\\
	1 & \text{if}\ x > 1,
	\end{cases},\qquad G(x) =
	\begin{cases}
	0 & \text{if}\ x < 0\\
	\frac{2}{3} & \text{if}\ 0 \le x < \frac{1}{4}\\
	1 & \text{if}\ x \ge \frac{1}{4}.
	\end{cases} \ui \ui     \hfill(1.7)
	\end{equation*}
	Measure $\nu$\ is discrete with atoms at $0\ \text{and}\ \frac{1}{4}$\ with saltus values $\frac{2}{3}\ \text{and}\ \frac{1}{3}$\ respectively. Let us find $\pi(\mu,\ \nu)$. For $C$,\ a closed subset of $R$,  $\nu(C) = \nu(C\cap \{0,\ \frac{1}{4}\})$. Recall $\pi(\mu,\ \nu) = \inf\{\varepsilon\ :\ \nu(C) \le \mu(C^{\varepsilon}) + \varepsilon\}$. It is sufficient to consider only the following three closed sets $C$ \ with {$\nu(C) > 0$}, namely, $C_1=\{0\} $, $\ C_2=\{1/4\}$  \text{and}$\ C_3=\{0, 1/4\}$.  We note $\nu(C_1) = \frac{2}{3}$. The least value of $\varepsilon$\ for which $\frac{2}{3} \le \mu(C_1^{\varepsilon}) + \varepsilon$\ is $\frac{1}{3}$. The least value of $\varepsilon$\ for which $\frac{1}{3} \le \mu(C_2^{\varepsilon}) + \varepsilon$\ is $\frac{1}{9}$. Let us examine $\{0,\ \frac{1}{4}\}^{\varepsilon}$. 
	It has to be of the form $(- \varepsilon,\ \varepsilon) \cup (\frac{1}{4} - \varepsilon,\ \frac{1}{4} + \varepsilon)$. If $\varepsilon \le \frac{1}{8}$, then $1 \leq \mu(C_3^{\varepsilon}) + \varepsilon = 4\varepsilon  $. Hence necessarily, $\varepsilon > \frac{1}{8}$. And in that case  $\mu(C_3^{\varepsilon}) + \varepsilon = \frac{1}{4} + 2\varepsilon$. This will be equal to $1$\ if $\varepsilon = \frac{3}{8}$. It follows $\pi(\mu,\ \nu) = \frac{3}{8}$.\\
	Furthe in this case $\ell(F, G)=3/8$.\\
\textbf{Remarks.}\\
1. The observation that  $\pi(\mu,\ \nu) =3/8$ and $\ell(F, G)=3/8$ leads to the question if they are equal always. $\ell(F, G) \leq \pi(\mu,\ \nu)$ is proved in Huber [4] (See Eq. (4.13) on page 34). We shall now strengthen this by the following result.\\
2. It is instructive to note that, if $\rho$\ is the uniform metric on the space of distribution functions, then $\rho(F,\ G) = \max\limits_{- \infty < x < \infty}|F(x) - G(x)| = \frac{3}{4}$. This raises the question if Prohorov distance gives the least distance ``in probability" between random variables distributed according to $F, \;G$ (measures $\mu,\;\nu$) . This indeed is true. We refer to Strassen [11] and Dudley [2]. \\\\
	\noindent\textbf{Theorem 1.4.}\\
Let $\mu,\ \nu$\ be two probability measures on \textbf{R} and let $F,\ G$\ be the corresponding distribution functions. Then $\ell(F,\ G) = \pi(\mu,\ \nu)$.\\ Proof.
We use definitions of $\ell(F,\ G)$ and $\pi(\mu,\ \nu)$.
Let $Q_1 = \{\varepsilon > 0 \big| \mu(A) \le \nu(A^{\varepsilon}) + \varepsilon\ \text{and}\;
\nu(A) \le \mu(A^{\varepsilon}) + \varepsilon\\\hspace*{5cm} \text{for all closed sets}\ A\ \text{ of the type}\ (- \infty,\ x] \}.$
Let $Q_2 = \{\varepsilon > 0 \big| \mu(A) \le \nu(A^{\varepsilon}) + \varepsilon\ \text{and}\;
\nu(A) \le \mu(A^{\varepsilon}) + \varepsilon\\\hspace*{5cm} \text{for all closed sets }\ A \}.$
We note $Q_2 \subset Q_1$. Hence $\ell(F,\ G) \le \pi(\mu,\ \nu)$.\hfill(1.8)
Suppose $\varepsilon \in Q_1$. This implies that for every $a < b$\ arbitrary, and every $k \ge 1$ \\\hspace*{.5cm} $F(b) - F(a - \frac{1}{k}) \le G(b + \varepsilon) + \varepsilon - \{G(a - \frac{1}{k} - \varepsilon) - \varepsilon\} \\\hspace*{3.5cm}\le G(b + \varepsilon) - G(a - \frac{1}{k} - \varepsilon) + 2\varepsilon$.\\
i.e.,\hspace*{3cm} $\mu((a - \frac{1}{k},\ b]) \le \nu((a - \frac{1}{k},\ b ]^{\varepsilon}) + 2\varepsilon.$
Since this is true for all $k$, we get, letting $k \rightarrow \infty$,\\\hspace*{3.5cm} $\mu([a,\ b]) \le \nu([a,\ b]^ {\varepsilon}) + 2\varepsilon$ \\
and similarly,\hspace*{1.3cm} $\nu([a,\ b]) \le \mu([a,\ b]^{\varepsilon}) + 2\varepsilon$ \\
both holding for all $a < b$\ and all $\varepsilon \in Q_1$.
If the closed set $A$\ is the complement of the union of $n$\ disjoint open intervals, then $A$\ would be the union of $2^n$\ disjoint closed intervals. Suppose $A$\ is the union of $m$\ disjoint closed intervals, say, $I_i,\ i = 1,\ 2,\ ..,\ m$. An $\eta > 0$\ can be found (which can be taken to be less than $\frac{\epsilon}{2m}$\ with no loss of generality) such that\\\hspace*{.5cm}$\mu(A) = \sum\limits_{i = 1}^m \mu(I_i) \le \sum\limits_{i = 1}^m \big(\nu(I_i^{\eta}) + 2\eta\big) \le \nu(A^{\eta}) + \varepsilon \le \nu(A^{\varepsilon}) + \varepsilon.$\\
Similarly,\hspace*{2cm} $\nu(A) \le \mu(A^{\varepsilon}) + \varepsilon.$\\\hspace*{.5cm}
These relations are true for all $\varepsilon \in Q_1$\ and for all closed sets which are the unions of finite number of disjoint closed intervals.\\\hspace{.5cm} Let now $A \subset \textbf{R}$\ be an arbitrary closed set. Hence $A = \bigcap\limits_{n = 1}^{\infty}J_n$\ where for $n = 1,\ 2,\ .....,\ J_n$\ is the union of $2^n$\ disjoint closed intervals and where $J_{n + 1} \subset J_n$.\\\hspace*{.5cm}$\mu(A) = \mu(\lim\limits_{n \rightarrow \infty}J_n) = \lim\limits_{n \rightarrow \infty}\mu(J_n) \le \lim\limits_{n \rightarrow \infty}\{\nu(J_n^{\varepsilon}) + \varepsilon\} \le \nu(A^{\varepsilon}) + \varepsilon.$\\ Similarly\hspace*{2.5cm}$\nu(A) \le \mu(A^{\varepsilon}) + \varepsilon.$\\\hspace*{.5cm}
These being true for all closed sets $A$, we conclude that $(\varepsilon \in Q_1) \Rightarrow (\varepsilon \in Q_2)$. Thus $Q_1 \subset Q_2$. Hence\\ $\pi(\mu,\ \nu) \le \ell(F,\ G).$\hfill(1.9)\\ 
The proof is completed by appealing to (1.8) and (1.9).\qed\\
\textbf{Remark.}\;Huber's proof of (1.8) is descriptive while our proof is constructive.The converse part may be new.\\
\textbf{Definition 1.6.}\\
A sequence $\{F_n\}$ of distribution functions is said to converge weakly to a distribution function $F$ if $F_n(x) \rightarrow F(x)$\ at all the continuity points $x$\ of $F$.\\
We shall now present  criteria for this weak convergence in terms of the Levy metric and the Prohorov metric.\\
\textbf{Theorem 1.5.}\\\hspace*{.5cm}
Let $F,\ F_n,\ n = 1,\ 2,\ ....$\ be distribution functions on the line. Then $\ell(F_n,\ F) \rightarrow 0$ \ if and only if $F_n(x) \rightarrow F(x)$\ at all the continuity points $x$\ of $F$.\\
Proof.\\\hspace*{.5cm}
Let  $\mathcal{C}_F$\ consists of all the continuity points of $F$ and let $u \in \mathcal{C}_F$.  If $\eta > 0$\ is given and if $\ell(F_n,\ F) \rightarrow 0$, then there exists $N$\ such that for all $\varepsilon < \eta$\ and $n \ge N$, 
$F(u - \varepsilon) - \varepsilon \le F_n(u) \le F(u + \varepsilon) + \varepsilon$. Take now limit as $n \rightarrow \infty$\ and get\\ $F(u - \varepsilon) -  \varepsilon \le \varliminf\limits_{n \rightarrow \infty}F_n(u) \le \varlimsup\limits_{n \rightarrow \infty}F_n(u) \le F(u + \varepsilon) + \varepsilon$. Now let $\varepsilon \rightarrow 0$\ and use fact that $u \in \mathcal{C}_F$\  to claim $\lim\limits_{n \rightarrow \infty}F_n(u) = F(u).$
Conversely, let $\varepsilon > 0$\ be given. Let $x \in \textbf{R}$\ be arbitrary. Given $\varepsilon$\ and $x$, we can find $u \in \mathcal{C}_F,\ u \in (x - \varepsilon,\ x + \varepsilon)$. This is possible since $\mathcal{C}_F^{\prime}$\ is at most countable. Let $\lim\limits_{n \rightarrow \infty}F_n(u) = F(u)$\ for every $u \in \mathcal{C}_F$\ . We have : there is $N$\ such that $|F_n(u) - F(u)| < \varepsilon,\ n \ge N$. Consequently, $F_n(x - \varepsilon) - \varepsilon \le F_n(u) - \varepsilon \le F(u) + \varepsilon - \varepsilon \le F_n(u) + \varepsilon \le F_n(x + \varepsilon) + \varepsilon$. This shows that $\ell(F_n,\ F) < \varepsilon$. Since $\varepsilon >0$\ is arbitrary, we conclude $\lim\limits_{n \rightarrow \infty}\ell(F_n,\ F) = 0.$\qed\\
It is now clear that the following holds.\\
\textbf{Theorem 1.6.}\\
Let $\mu,\ \mu_n \in \mathcal{M},\ n \ge 1$. Then $\  (\pi(\mu_n,\ \mu)\ \rightarrow \ 0) \Rightarrow \ (\mu_n\ \ ^{\underrightarrow{w}}\ \ \mu)$.\\
\textbf{Remark.}\\
	If $\mu_n$  and  $\mu$ are probability measures on \textbf{R} determined by distribution functions $F_n$ and $F$ then by Theorems 1.3 - 1.5 we have the result  $(\pi(\mu_n,\ \mu)\ \rightarrow \ 0) \Leftrightarrow (\ell(F_n,\ F)\ \rightarrow \ 0) \Leftrightarrow (\mu_n\ \ ^{\underrightarrow{w}}\ \ \mu)$. Refer to Theorem 2.2 below.
\section{Tightness and weak compactness.}
In this Section we define and investigate the tightness and weak compactness of probability measures defined on the Borel $\sigma$-field of a topological space $M$. We state and prove Prohorov's theorem.\\
\textbf{Definition  2.1.}\\
A family $(\mu_{\alpha})$\ of probability measures  on the Borel $\sigma$-field $\textit{\Large{m}}$\ (i.e., the $\sigma$-field generated by the open subsets) of a topological space $\textbf{M}$\ is said to be tight if, given $\varepsilon > 0$, a compact set $K_{\varepsilon}$\ can be found such that $\mu_{\alpha}(K_{\varepsilon}) > 1 - \varepsilon$\ for every $\alpha$\ in the index set.\\
\noindent\textbf{Theorem 2.1}.\\\hspace*{.5cm}
Every probability measure $\mu$\ on the Borel $\sigma$-field $\textit{\Large{m}}$\ of a complete and separable metric space \textbf{M}\ is tight.\\
\\\hspace*{.5cm}Let $(a_1,\ a_2,\ ...)$\ be a separability set for $(\textbf{M},\ d)$. Denote by $\overline{S}_{n,j}$\ the closed sphere with center at $a_j$\ and radius $\frac{1}{n}$. 
Given $\varepsilon > 0$, we can find $k_n$\ such that $\mu(B_n) > 1 - \frac{\varepsilon}{2^{n + 1}}$\ where $B_n = \bigcup\limits_{j = 1}^{k_n}\overline{S}_{n,j}$. This is possible since, for each $n,\ \textbf{M} = \bigcup\limits_{j = 1}^{\infty}\overline{S}_{n,j}$. Define $K = \bigcap\limits_{n = 1}^{\infty}B_n$. Since \textbf{M}\ is a complete metric space, it follows $K$\ is a compact set . The tightness of $\mu$\ is now immediate since $\mu(K^{\prime}) \le \sum\limits_{n = 1}^{\infty}\frac{\varepsilon}{2^{n + 1}} < \varepsilon$.\qed\\
Second proof.\\\hspace*{.5cm}
	Let $\mathscr{Y}$\ be the Borel $\sigma$-field of \textbf{Y} $= [0,\ 1]^{\infty}$. \\\hspace*{.5cm} We know there exists a homeomorphic map $\varphi$, mapping $\textbf{M}$\ on to a Borel subset of $\textbf{Y}$.
	Appealing to the definition of a tight measure, we note that every probability measure on the Borel $\sigma$-field of a compact metric space is tight. (Thus every probability measure on $\mathscr{Y}$\ is tight, since, by Tychonoff theorem, \textbf{Y}\ is a compact set.) 
	We therefore assume $\textbf{M}$\ is not compact, as otherwise there is nothing to prove.\\\hspace*{.5cm}This assumption implies that the Borel set $\varphi(\textbf{M})$\ can not be a closed set. For, were it so, being a closed subset of the compact set \textbf{Y}, $\varphi(\textbf{M})$\ would be a compact set and that would imply $\textbf{M}$\ is a compact set. The assumption implies also that $\varphi(\textbf{M})$\ is a proper subset of \textbf{Y}.
	\\\hspace*{.6cm}The Borel $\sigma$-field of $\varphi(\textbf{M})$\ endowed with its relative topology would be $\mathscr{Y} \cap \varphi(\textbf{M})$. \\\hspace*{.5cm} The Borel $\sigma$- field generated by the relative topology of $\varphi(\textit{M})$\ is $\mathscr{Y} \cap \varphi( \textbf{M})$. It is also equal to $\mathscr{Y} \cap \varphi( \textit{\Large{m}})$. Thus if $D \in \mathscr{Y}$\ then there is $E \in \textit{\Large{m}}$\ such that $D \cap \varphi(\textbf{M}) = \varphi(E)$. \\\hspace*{.5cm}Define measure $\nu$\ on $\varphi(\textit{\Large{m}})$: $\nu(\varphi(E)) = \mu(E)$\ for every $E \in \textit{ \Large{m}}$. Define measure $\nu^*$\ on $\mathscr{Y}$: if $D \in \mathscr{Y}$, write $\nu^*(D) = \nu(D \cap \varphi(\textbf{M})) = \mu(E)$. 
	\\\hspace*{.5cm}We note that each member of the collection of sets\\\hspace*{1.6cm} $\mathcal{C} = \{\varphi(C) \big| C \subset \textbf{M},\ C\ \text{compact}\}$\ is a compact subset of $\varphi(\textbf{M})$. We note that every compact subset of $\varphi(\textbf{M})$\ is a compact subset of \textbf{Y}, by reason of $\varphi(\textbf{M})$ having the inherited metric.\\\hspace*{.5cm}Are there any other compact subset of $\varphi(\textbf{M})$? No. For, if $E \subset \varphi(\textbf{M})$\ is compact, then $\varphi^{-1}(E)$\ would be a compact subset of $\textbf{M}$. Hence $E \in \mathcal{C}$.\\\hspace*{.5cm}
	As defined above, $\nu^*$\ is a tight measure. Hence, given $\varepsilon > 0$, there is a compact set $K$\ with $\nu^*(K) > 1 - \varepsilon$. But all relevant compact sets are in $\mathcal{C}$. Thus $K \in \mathcal{C}$. This implies that $D = \varphi^{-1}(K)$\ is a compact subset of $\textbf{M}$ and $\mu(D) > 1 - \varepsilon$. The proof that $\mu$\ is tight is now complete.\qed \\\\
The following is a converse to the result in Theorem 1.6 when the metric space is separable.\\
	\noindent\textbf{Theorem 2.2.}\\
Let $(M,\ d)$\ be a separable metric space. If {probability measures}\; $\mu,\ \mu_n \in \mathcal{M}$\ and if $\mu_n\ \  ^{\underrightarrow{w}}\ \ \mu$, then $\pi(\mu_n,\ \mu)\ \rightarrow\ 0$.\\
Proof.
\hspace*{.6cm} Step 1 \\
\hspace*{.5cm} Fix $\varepsilon > 0$. The theorem will stand proved if we can find $N$\ such that for all $n \ge N,\ \pi(\mu_n,\ \mu) < \varepsilon$. This will follow if we show\\\hspace*{2cm}  $\mu(B) \le \mu_n (B^{\varepsilon}) + \varepsilon$\ $\&\ \mu_n(B) \le \mu(B^{\varepsilon}) + \varepsilon$\ for all $n \ge N$\\\hspace*{7.5cm} and all Borel sets $B$.\hfill(2.1)\\\hspace*{.5cm}Let $S = \{a_{j}\}$\ be a separability set for $(\textbf{M},\ d)$. Let $\delta > 0, \delta < \frac{\varepsilon}{3}$. As argued in Theorem 1.2  we can find, for each $j$, a closed sphere $\overline{S}_{j}\ $ with center at $a_{j}$\ and radius less than $\frac{1}{2}\delta$\ and such that $\mu(\partial \overline{S}_{j}) = 0$,
Since $S$\ is dense in \textbf{M},  $\textbf{M} = \bigcup \limits_{j = 1}^{\infty}\overline{S}_{j}$. Determine $k$\ such that\\\hspace*{2.5cm} $\mu(A) > 1 - \delta$\ where $A = \bigcup\limits_{j = 1}^k\overline{S}_j$.\hfill(2.2)\\\hspace*{.5cm}We note $\mu(A^{\prime}) < \delta.$ \\\hspace*{.5cm}For the Borel set $A$\ in (2.2), $\partial A \subset \bigcup\limits_ {j = 1}^k\partial \overline{S}_j$. Hence $\mu(\partial A) = 0$. Since $\partial A^{\prime} = \partial A$, $\mu(\partial A^{\prime}) = 0.$\\
\hspace*{.5cm}Let now $B$\ be an arbitrary Borel set. \\
Case 1. $B \cap A = \emptyset$.\ So $B \subseteq A^{\prime}$. Since $\mu(\partial A^{\prime}) = 0$, there exists $N_1$\ such that $|\mu_n(A^{\prime}) - \mu(A^{\prime})| < \delta\ \text{for all}\ n \ge N_1$. In this case (2.1) holds since\\\hspace*{2cm}
$\mu(B) \le \mu(A^{\prime}) < \delta <  \mu_n(B^{\varepsilon}) + \varepsilon$\ and \\ \hspace*{2cm}$\mu_n(B) \le \mu_n(A^{\prime}) \le \mu(A^{\prime}) + \delta \le 2\delta \le \mu(B^{\varepsilon}) + \varepsilon.$\\ 
Case 2. $B \cap A \neq \emptyset$. Let $J \subset \{1,\ 2,\ \ldots,\ k\}$\ be such that $j \in J$\ if and only if $B \cap \overline{S}_j \neq \emptyset$. Let $E = \bigcup \limits_{j \in J}\overline{S}_j$\ and note $\mu(\partial E) = 0$.\; So there is $N_2$\ such tat for all $n \ge N_2,\ui\ui |\mu_n(E) - \mu(E)| < \delta$.\\ Choose $\delta$\ to satisfy the further condition :\  $\mu(\partial B^{\delta}) = 0$. So there exists $N_3$\ such that, for all $n \ge N_3$,\\ \hspace*{4cm}$|\mu_n(B^{\delta}) - \mu(B^{\delta})| < \delta.$ \hfill(2.3)\\ We note, if $j \in J$, $B^{\delta} \supset B \cap \overline{S} _j$ and hence  $B \cap E \subset B^{\delta}$.  \\\hspace*{.5cm}Now, \\\hspace*{2cm}$\mu(B) = \mu(B \cap A) + \mu(B \cap A^{\prime}) \le \mu(B \cap E) + \mu(A^{\prime}) \le \mu(B^{\delta}) + \delta\\ \hspace*{2.9cm} \le \mu_n(B^{\delta}) + 2\delta \le \mu_n(B^{\varepsilon}) + \varepsilon.$\\ Again,\\ \hspace*{2cm}$\mu_n(B) \le \mu_n(B^{\delta}) \le \mu(B^{\delta}) + \delta \le  \mu(B^{\varepsilon}) + \varepsilon.$\\\hspace*{.5cm}With this the proof of the theorem is complete.\qed\\
$\textbf{Y} = [0,\ 1]^{\infty}$\ is  endowed with the product topology (equivalent to the topology induced by the metric $\rho$). $(\textbf{Y},\ \rho)$\ is a compact metric space, complete and separable.  Since $\textbf{Y}$\ is compact, every family of probability measures on the Borel $\sigma$-field $\mathscr{Y}$\ of \textbf{Y}\ is tight.\\ 
Define projection operators $\wp_{j_1,\ j_2,\ ...,\ j_k}(\textbf{y}) = (y_{j_1},\ y_{j_2},\ ...,\ y_{j_k})$. Since convergence in the $\rho$-metric is co-ordinatewise convergence, it is clear that all projection operators are continuous. \\
We now discuss the sequential compactness of probability measures.\\
\textbf{Theorem 2.3}.\\\hspace*{.5cm}Every infinite sequence of probability measures on $\mathscr{Y}$\ contains a weakly convergent subsequence.\\ Proof.\\\hspace*{,5cm}
Let $\mu_n \in A$\ be a tight sequence. Given $\varepsilon > 0$, there exists then a compact set $K \subset \textbf{Y}$\ such that $\mu_n(K) > 1 - \varepsilon$\ for all $n \ge 1$. Define $C_r = \wp_{1,\ 2,\ ...,\ r}(K)$. Note that it is a compact subset of $\textbf{R}^k$\ and that $\wp_{1,\ 2,\ ...,\ r}^{-1}C_r \supset K$. Hence $\mu_n(\wp_{1,\ 2,\ ...,\ r}^{-1}C_r) > 1 - \varepsilon$. Thus  $\mu_n(\wp_{1,\ 2,\ ...,\ r}^{-1}), n \ge 1$\ is a tight sequence of probability measures on $\textbf{R}^k$. Hence it contains a weakly convergent subsequence (ref. {p. 85,} Pakshirajan [7]).\\\hspace*{.5cm}For $r = 1$, denote the weakly convergent  subsequence by $\mu_{1,n} \wp_1^{-1}$. Now, ($\mu_{1,n}$), being a subsequence of $(\mu_n)$, is tight. This implies the tightness of $(\mu_{1,n}\wp_{1,2}^{-1})$, which is a sequence of probability measures on $\textbf{R}^2$. Hence it contains a weakly convergent subsequence. Denote it by $\mu_{2,n}\wp_{1,2}^{-1}$. Arguing on these lines, we arrive at $\mu_{3,n},\ \mu_{4,n},\ ...$. We note, that for each $j$, sequence $(\mu_{j + 1,n})$\ is a subsequence of $(\mu_{j,n})$. The diagonal sequence $(\mu_{n,n})$\ has the property that, (i) it is tight , since it is a subsequence of the $\mu_n$-sequence and (ii) for every $k \ge 1$, $(\mu_{n,n}\wp_{1,\ 2,\ ...,\ k}^{-1})$\ is weakly convergent.\\\hspace*{.5cm}Let $\mu_{n,n}\wp_{1,\ 2,\ ...,\ k}^{-1}\  ^{\underrightarrow{w}}\ \alpha_k$. If $\mathscr{R}^k$\ denotes the Borel $\sigma$-field of $\textbf{R}^k$, then $\mathscr{Y}_k = \wp_{1,\ 2,\ ...,\ k}^{-1}(\mathscr{R}^k)$\ would be a sub $\sigma$-field of $\mathscr{Y}$. We note $\wp_{1,\ 2,\ ....\ k}^{-1}(E) = \wp_{1,\ 2,\ ...,\ k +1}^ {-1}(E \times \textbf{R})$. This shows that $\mathscr{Y}_k \subset \mathscr{Y}_{k + 1}$.
Let $\lambda_{k + 1}$\ be the projection operator mapping $\textbf{R}^{k + 1}$\ on to $\textbf{R}^k$. Then $\lambda_{k + 1}\wp_{1,\ 2,\ ...,\ k + 1} = \wp_{1,\ 2,\ ...,\ k}$. Hence ($\mu_{n,n}\wp_{1,\ 2,\ ....,\ k}^{-1}\ ^{\underrightarrow{w}}\ \alpha_k$) $\Leftrightarrow$\ $(\mu_{n,n}\wp_{1,\ 2,\ ....,\ k + 1}^{-1}\lambda_{k + 1}^{-1}\ ^{\underrightarrow{w}}\ \alpha_k)$. But\\ $\mu_{n,n}\wp_{1,\ 2,\ ....,\ k + 1}^{-1}\ ^{\underrightarrow{w}}\ \alpha_{k + 1}$. Hence $\alpha_k = \alpha_{k + 1}\lambda_{k + 1}^{-1}$. This shows that the $\alpha_k$s defined on the $\sigma$-fields $\mathscr{Y}_k$\ form a consistent family of measures. Hence (ref. Note (under Theorem 1.11.1, p. 51, Pakshirajan [7]) there exists a probability measure $\mu$\ on $\mathscr{Y} (= \sigma(\bigcup_{k = 1}^{\infty}\mathscr{Y}_k))$\ such that $\mu\wp_{1,\ 2,\ ...,\ k}^{-1} =  \alpha_k,\ k \ge 1$.\\\hspace*{.5cm}That $\mu_{n,n}\ ^{\underrightarrow{w}}\ \mu$\ follows now from the following Lemma with $\nu_n=\mu_{n, n} $ and $ \nu=\mu$.\qed\\
The following Lemma  provides a criterion for the weak convergence of a family of probability measures to a probability measure.\\
\textbf{Lemma 2.1.}\; Let ($\nu_n$) be a sequence of probability measures defined on the Borel $\sigma$-field $\mathscr{Y}$\ defined above. Let ($\nu_n$) satisfy the following three conditions:\\
(i) ($\nu_n$) is tight\\
(ii) for every $k \geq 1$ the sequence ($ \nu_n\;\pi_{1, 2, \ldots, k}^{-1}$) has a weak limit where $\pi_{1,\ 2,\ldots, k} $ is a projection operator on compact subsets of $Y$  to $R^k$, and\\
(iii) there exists a probability measure $\nu$ having the weak limits in (ii) for its finite dimensional distributions.\\
Then $\nu_n\ ^{\underrightarrow{w}}\ \nu$.\\
\textbf{Proof} \;Let $K \subset \textbf{Y}$\ be compact  and let $C_k = \pi_{1,\ 2,\ ...,\ k}(K)$ which  is  now compact. Consider the sequence of measures ($ \nu_n\;\pi_{1,\ 2,\ ...,\ k}^{-1}$). Note that $\nu_n(K) \leq \nu_n\;\pi_{1,\ 2,\ ...,\ k}^{-1} \;(C_k) $ and hence\\
\hspace*{4cm}$\varlimsup\limits_{n \rightarrow \infty} \nu_n(K)\leq \varlimsup\limits_{n \rightarrow \infty} \nu_n\;\pi_{1,\ 2,\ ...,\ k}^{-1}(C_k).$ \;\hfill (2.4)
\\We extend the result in (2.4) to  an arbitrary closed set $C$. Given $\varepsilon > 0$ by (i) we can find a compact set $K \subset \textbf{Y}$ such that, for all $n \geq 1$, 
$\nu_n(K) > 1- \varepsilon$. Now\\
$\nu_n(C) =  \nu_n(C\cap K) + \nu_n(C\cap K^\prime) \leq  \nu_n(C\cap K) + \nu_n( K^\prime) \leq  \nu_n(C\cap K) + \varepsilon. $ \\
Then observing that $C\cap K$ is a compact set and using (ii), (iii) and Theorem 1.1, we have from (2.4)
$$ \varlimsup_{n \longrightarrow\infty} \nu_n(C) \leq \varlimsup_{n \longrightarrow\infty} \nu_n(C\cap K) + \varepsilon \leq \nu(C\cap K)+\varepsilon \leq \nu(C) + \varepsilon.$$
This being true for all $\varepsilon > 0$ we get $\varlimsup_{n \longrightarrow\infty} \nu_n(C) \leq \nu (C)$ for every closed set $C$. Then by Theorem 1.1,  $\nu_n \stackrel{w}{\longrightarrow} \nu.$\\
\textbf{Remark.} The converse to the  result in the above Lemma 2.1 1 is also true.\\
 Next we discuss the separability of the space $(\mathcal{M},\ \pi)$. \\
\textbf{Theorem 2.4}.\\\hspace*{.5cm}If $(\textbf{M},\ d)$\ is separable then so is $(\mathcal{M},\ \pi)$.\\ Proof.\\\hspace*{.5cm}Let $\mu \in \mathcal{M}$\ be arbitrary. Let $(a_n)$\ be a separability set for $\textbf{M}$. For each $n = 1,\ 2,\ ...$, the closed spheres $S_j = S(a_j,\ n)$ with center at $a_j$\ and $\frac{1}{2n} < radius  < \frac{1}{n};\ j = 1,\ 2,\ ...$\ is a cover for \textbf{M}. 
Define $ V_{n,j} = V_j: V_1 = S_1;\ V_2 =  S_1^{\prime} \cap S_2;\ V_3 = S_1^{\prime} \cap S_2^{\prime} \cap S_3$\ and so on. We note that the diameter of each $V_{n,j}$\ is $\le \frac{1}{n}$, that for each $n$, the sets $V_{n,j},\ j \ge 1$\ is a disjoint collection and that $\bigcup\limits_{j = 1}^{\infty}V_{n,j} = \textbf{M}$. 
Note each $V_{n,j}$ contains an open set. Let $b_{n,j} \in V_{n,j}$\ be chosen from the separability set  and fixed.  Define discrete probability measure $\mu_n\ :\ \mu_n(\{b_{n,j}\}) = \mu(V_{n,j})$.  To claim $\mu_n\ ^{\underrightarrow{w}}\ \mu$, we show that condition (iv) of Theorem 1.1 is satisfied. Let $f$\ be uniformly continuous. Given $\varepsilon > 0$, we can find $N$\ such that $|f(x) - f(y)| < \varepsilon$\ whenever $d(x,\ y) < \frac{1}{N}$. This is possible since $f$\ is uniformly continuous.\\\hspace*{.5cm}For $x \in V_{j,n}, d(x,\ b_{n,j}) \le \frac{1}{n}$. Hence for all $n > N$,\\
$\big|\int\limits_{\textbf{M}}f(x)\ud\mu_n(x) - \int\limits_{\textbf{M}}f(x)\ud\mu(x)\big| = \big|\sum\limits_{j = 1}^{\infty}\ \int\limits_{V_{n,j}}\{f(x) - f(b_{n,j})\}\ud\mu(x)\big|\\\hspace*{5.4cm}\le \sum\limits_{j = 1}^{\infty}\int\limits_{V_{n,j}}|f(x) - f(b_{n,j})|\ud\mu(x) \\\hspace*{5.4cm} < \varepsilon \sum\limits_{j = 1}^{\infty}\int\limits_{V_{n,j}}\ud\mu(x) < \varepsilon$. \\ This shows that $\int\limits_{\textbf{M}}f(x)\ud\mu_n(x) \rightarrow \int\limits_{\textbf{M}}f(x)\ud\mu(x)$. In other words (ref. Theorem 1.1) we have shown $\mu_n \ ^{\underrightarrow{w}} \ \mu$. i.e., we have shown that every member of $\mathcal{M}$\ is the weak limit of a sequence of measures, whose supports are sets with a countable number of members from the separability set. In turn these measures are the weak limits of measures concentrated on a finite number of points. To summarise, every $\mu \in \mathcal{M}$\ is the weak limit of a sequence of measures with support in a finite   subset of the separability set. Hence $(\mathcal{M},\ \pi)$\ is a separable metric space.\qed\\
\\Let $\textit{\Large{m}}$\ be the Borel $\sigma$-field \ of a complete and separable metric space $(\textbf{M},\ d)$. Let $(\mathcal{M},\ \pi)$\ be the metric space of all the probability measures on
$\textit{\Large{m}}$, $\pi$\ being the Prohorov metric. \\
\textbf{Definition 2.2.}\\\hspace*{.5cm}A family $\mathcal{F}$\ of probability measures on $\textit{\Large{m}}$\ is said to be relatively sequentially compact if every sequence in it contains a weakly convergent subsequence.\\\hspace*{.5cm}Note.\\\hspace*{.5cm}Saying that $\mathcal{F} \subset \mathcal{M}$\ is relatively sequentially compact is equivalent to saying tha $\overline{\mathcal{F}}$\ is compact.\\
We shall now present the main result of this Section.\\
\textbf{Theorem 2.5.} (Prohorov [10]).\\ A family $\mathcal{F} \subset \mathcal{M} $ is tight if and only if its closure in $(\mathcal{M},\ \pi)$\ is compact.
\\ Proof.\\ \hspace*{.5cm}  Let $\varepsilon > 0$\ be given. Tightness of the sequence ($\mu_n$) implies that there exists a compact set $K = K_{\varepsilon}$\ such that $\mu_n(K) > \varepsilon$\ for all $n$. 
There then exists a continuous function $\varphi$\ mapping $\textbf{M}$\ on to a Borel subset of $\textbf{Y}$\ such that $\varphi$\ is one-to-one\ and $\varphi^{-1}$\ defined on $\varphi(\textbf{M})$\ is continuous. Note that $\varphi(K)$\ is a compact subset of $\textbf{Y}$. Define probability measures $\nu_n\ :\ \nu_n = \mu_n\varphi^{-1} $. $(\nu_n)$\ is a tight sequence of measures on $(\varphi(\textbf{M}),\ \rho)$\ since $\nu_n (\varphi(K)) = \mu_n\varphi^{-1}(\varphi(K)) = \mu_n(K) > 1 - \varepsilon$. 
The $\nu_n$s can be thought of as defined on $\mathscr{Y}$\ in a natural way. By  Theorem 2.3 the $\nu_n$-sequence contains a weakly convergent sub sequence, say, $(\nu_{n_k})$ converging to, say, $\nu$. Hence\\ $1 - \varepsilon < \varlimsup\limits_{k \rightarrow \infty} \nu_{n_k}(\varphi(K_{\varepsilon})) \le \nu(\varphi(K_{\varepsilon}))$. This being true for every $\varepsilon > 0$, we get $\nu(\varphi (\textbf{M})) = 1$. Define $\mu = \nu\varphi$. If $D \subset \textbf{M}$\ is an open set, then $\varphi(D)$\ is an open subset of $\varphi(\textbf{M})$ and $\varliminf\limits_{k \rightarrow \infty}\mu_{n_k}(D) = \varliminf\limits_{k \rightarrow \infty}\nu_ {n_k}(\varphi(D)) \ge \nu(\varphi(D)) = \mu(D)$. This shows $\mu_{n_k}\ ^{\underrightarrow{w}}\ \mu$.\\\hspace*{.5cm}  Let $\mathcal{M}$\ denote the totality of all probability measures on $\textit{\Large{m}}$. Let $\mathcal{F} \subset \mathcal{M}$\ be relatively sequentially compact. Let $(a_n)$\ be a separability set for \textbf{M}. Let $S(a_n,\ \delta)$\ be the open sphere with center at $a_n,\ n = 1,\ 2,\ ...$\ and radius $\delta$. Clearly $\textbf{M} = \bigcup\limits_{n = 1}^{\infty}S(a_n,\ \delta)$. Hence for each $\mu \in \mathcal{M}$\ we can find integer $q(\mu,\ \delta)$\ such that $\mu(\bigcup\limits_{j = 1}^{q(\mu,\ \delta)}\overline{S(a_j,\ \delta)}  \ ) > 1 - \delta$. We claim there exists $q(\mathcal{F},\ \delta)$\ such that \hspace*{1cm} $\mu(\bigcup\limits_{j = 1}^{q(\mathcal{F},\ \delta)} \overline{S(a_j,\ \delta)}\ ) > 1 - \delta$\ for all $\mu \in \mathcal{F}.$\hfill(2.5) \\\hspace*{.5cm}If this claim is not admitted, then whatever number $r$\ we choose, \\$\mu(\bigcup\limits_{j = 1}^r\overline{S(a_j,\ \delta)}\ ) \le 1 - \delta$\ for infinitely many $\mu \in \mathcal{F}$. Let $(\mu_n)$\ be a sequence for which this inequality holds for each $n$. Since the family is relatively sequentially compact, sequence $(\mu_n)$\ contains a weakly convergent subsequence, say, $(\mu_{n_k})$\ converging to, say, $\mu$. Since the union set is a closed one and since $\mu_{n_k}\ \ ^{\underrightarrow{w}}\ \ \mu$\ and consequently $\varlimsup\limits_{k \rightarrow \infty}\mu_{n_k}(\bigcup\limits_{j = 1}^r\overline{S(a_j,\ \delta)}\ ) \le \mu(\bigcup\limits_{j = 1}^r\overline{S(a_j,\ \delta)}\ )$ we get :\  $\mu(\bigcup\limits_{j = 1}^r\overline{S(a_j,\ \delta)}\ ) \le 1 - \delta$, leading to $\mu( \textbf{M}) \le 1 - \delta$\ which is absurd. Hence (2.5) holds.\\\hspace*{.5cm}Define $K = \bigcap \limits_{n = 1}^{\infty}\bigcup\limits_{j = 1}^{q(\mathcal{F},\ \frac{\varepsilon}{2^n})} \overline{S(a_j,\ \frac{\varepsilon}{2^n})}$\ and it can be proved that $K$\ is a compact set. Since (2.5) holds for every $\mu \in \mathcal{F}$, 
$\mu(K^{\prime}) \le \sum\limits_{n = 1}^{\infty}\frac{\varepsilon}{2^n} < \varepsilon$. Thus $\mu(K) > 1 - \varepsilon$\ for every $\mu \in \mathcal{F}$. The proof is now complete that the family $\mathcal{F}$\ is tight.\qed\\
Prohorov theorem establishes equivalence of tightness of  $\mathcal{D} \subset \mathcal{M}$ and  compactness of its closure in  ($\mathcal{M},\; \pi$). 
However in specific examples  of complete and separable spaces  it is difficult to say which of the above  equivalent  properties is easy to check, although tightness appears to be easier to check than the other. We now give an example  where this is indeed true.
Consider  $\textbf{M}=\textbf{C}[0,\ 1]$, the space of all real valued continuous functions $f$ on [0, 1] equipped with the uniform metric $\rho$  such that $f(0)=0$. The metric space ($\textbf{M}, \rho$) is complete and separable.
Denote the Borel $\sigma$-field on this metric space by $\textit{\Large{m}}$.\\
Let $\mathcal{D}$ be a compact subset of $\textbf{M}$ and let $\mu$ be a measure on  $\textit{\Large{m}}$. For Borel set $A \in \textit{\Large{m}}$ and function $f \in \mathcal{D}$ define the measure $\nu_f(A)=\mu(A - f)$. Note that $\nu_f$ is a probability  measure and the collection $\mathcal{E}=\left\lbrace \nu_f\right\rbrace \subset \mathcal{M}$. We shall show that $\mathcal{E}$ is tight: i.e., given $\varepsilon > 0$ there exists a  compact subset $K_\varepsilon^* \subset \textbf{M}$ such that $\nu_f (K_\varepsilon^*) > 1 - \varepsilon$ for every $\nu_f \in \mathcal{E}.$\\
In view of Theorem 2.1, given an $\varepsilon > 0$, there exists a compact subset $K_\varepsilon$ of \textbf{M} such that $\mu(K_\varepsilon) > 1- \varepsilon$. Without loss of generality we may assume that the element  $x:\;x(t)=0$ lies in $K_\varepsilon$.\\
Now define set $K_\varepsilon^*\;=\;K_\varepsilon \oplus \mathcal{D}=\{f + g|\;f \in \mathcal{D} ;\; g \in K_\varepsilon\}$. Note that this is a compact subset of $\textbf{M}$ and that $ K_\varepsilon \cup \mathcal{D} \subset K_\varepsilon \oplus \mathcal{D}\;=\;K_\varepsilon^*$.\\
Then for arbitrary  $f \in \mathcal{D}$, we observe 
$$\nu_f (K_\varepsilon \oplus \mathcal{D}) \geq  \nu_f( K_\varepsilon \oplus \{f\}) = \mu(K_\varepsilon \oplus \{f\}\ - f)= \mu(K_\varepsilon) > 1-\varepsilon.$$ \\ 
\textbf{Remarks.}\\ 1. Pakshirajan [6] proved the Prohorov theorem in $D[0, 1]$,  the space of real functions $x$ in $[0, 1]$ that are right-continuous and have left-hand limits. We shall discuss this in some detail in the following.\\
2. Pakshirajan [8] also proved Prohorov theorem in Banach spaces with Schauder bases.\\
3. We refer to Preiss [9] and the references therein, for  examples of spaces in which the Prohorov theorem is not valid.\\
Let $D[0, 1]$,  the space of real functions $x$ in $[0, 1]$ that are right-continuous and have left-hand limits. It follows that if $x, y \in D$
and if $x (t) \in  \{y(t), y(t-)\}$  for each $t$ belonging to a countable dense subset, then $x(t)= y(t$) for all $ t$. Given $t$, let
$t^*$ stand for $t$ or $t-$. For $0\leq t, t_1,t_2, \ldots, t_k \leq 1$,  define $\pi_{t_1^*, t_2^*, \ldots, t_k^*}\;x=(x(t_1^*), x(t_2^*), \ldots, x(t_k^*))$ mapping $D$ into $R^k$. Define operators mapping subsets of $D$ into subsets of the Euclidean
space of appropriate dimension: $\; \pi_{[t]}(\{x\})=\{x(t), x(t-)\},\; \pi_{[t_1, t_2. \ldots, t_k]}(\{x\})=\times_1^k\;\pi_{[t_i]}(\{x\})$.
For $A \subset D$ define $\pi_{[t_1, t_2. \ldots, t_k]}A=\bigcup_{x\in A}\;\pi_{[t_1, t_2. \ldots, t_k]}(\{x\})$\\
Let $T=\{t_n, n\geq 1\}$ be a fixed countable dense subset of [0, 1].  Let $\tau$ be a metric for $D$. Let $\mathcal{B}_\tau$ be the resulting  Borel $\sigma$-field of $D$.\\
\textbf{Definition 2.3.}\\ The metric $\tau$ is said to be regular if \\
(a) for every choice of $k\geq 1$  and every choice of $t$-values $\{t_1, t_2, \ldots, t_k\}$\\
(i) $\pi_{t_1^*, t_2^*, \ldots, t_k^*}$ are $\mathcal{B}_\tau$-measurable, and \\
(ii) for every compact set $K$,  $\pi_{[t_1, t_2. \ldots, t_k]} K$ is a closed subset of $R^k$, and \\
(b) $\tau(x_n, x)\rightarrow 0$ implies that the limit set of each of the two sequences $\{x_n(t)\}$  and $\{x_n(t-)\} $ is contained in the set $ \{x(t), x(t-)\}$ for every $t\in T.$\\
Then we have the following results. We assume that  $\tau$ is a regular metric. \\\\
\textbf{Theorem 2.6} (Pakshirajan [6]) If $K$ is a compact subset in $(D, \tau)$ then, for every $t\in T$, $A=\pi_{[t]}K$ is a compact subset of $ R$.\\
\textbf{ Remark.} An immediate consequence of the theorem is: $\pi_{[t_1, t_2, \ldots, t_k]\;}K $ is a compact subset of $R^k$ for every choice
of $k\geq 1$ and every choice of $\{t_1, t_2, \ldots, t_k\} \subset  T$.\\
\textbf{Theorem 2.7.}  (Pakshirajan [6]) Let $K$ be a compact subset of $D$. Let $E_n=\pi_{[t_1, t_2, \ldots, t_n]}K$ and $Q_n=\pi_{[t_1, t_2, \ldots, t_n]}^{-1}E_n$. Then $K=\bigcap_1^\infty\;Q_n.$\\\
In the following discussion all probability measures  $ \mu_n, \nu_n, m $ with or without suffix will be assumed to be defined on  $\mathcal{B}_{\tau}.$
By the finite dimensional distributions (fdd) of a probability measure $\mu$ we understand the family of probability measures induced on $\mathcal{R}^k$, the Borel $\sigma$-field of $R^k$, by $\pi_{[t_1, t_2, \ldots, t_k]}$ for every choice of $k\geq 1$ and every choice of $t_1, t_2, \ldots, t_k$.\\
Recalling the definition of tightness of probability measures, we have the following results.\\
\textbf{Theorem 2.8.} (Pakshirajan [6]) Two probability measures $\mu$ and $\nu$ with the same fdd are identical if one of them is tight.\\
\textbf{Theorem 2.9.} (Pakshirajan [6]) (a) If the fdd of $\mu_n$ converge weakly to the corresponding ones of $\mu$, then $\limsup_{n\rightarrow \infty}\mu_n\;(K)\leq \mu (K)$ for every compact set $K$.\\
(b) If $(\mu_n)$ is a tight sequence  then $\mu$ is a tight measure and $\mu_n$ converges weakly to $\mu$.\\
Let$\mathcal{B}^* \subset \mathcal{B_\tau}$  denote the minimal $\sigma$-field with respect to which the projections $\pi_t,  t \in T$  are measurable. Then by Theorem 2.7, all compact subsets of $(D; \tau))$ are in $\mathcal{B}^*$. We then have \\
\textbf{Theorem 2.10.} (Pakshirajan [6])  Let $(\mu_n)$ be a sequence of probability measures on $\mathcal{B}^* $ such that for every $k\geq 1$, the sequence $(\alpha_{n, k}=\mu_n \ \pi_{[t_1, t_2, \ldots, t_k]}^{-1}, n \geq 1)$  is weakly convergent. Denote the limit measure by $\alpha_k$.  Then there exists a unique probability measure $\mu$ on $\mathcal{B}^*$ such that
$\mu \ \pi_{[t_1, t_2, \ldots,t_k]}^{-1}=\alpha_k.$\\We now state and prove Prohorov theorem on $D[0, 1]$.\\
\textbf{Theorem 2.11.}  Every tight sequence $(\mu_n)$ of probability measures on $\mathcal{B}_\tau$ admits of a weakly convergent subsequence.\\
\textbf{Proof} Let $\varepsilon > 0$ be given. Since the sequence $(\mu_n)$ is tight, a compact set $K \subset D$ can be found such that
$\mu_n(K) >1-\varepsilon$ for all $n\geq 1$. Then $C_k=\pi_{[t_1, t_2, \ldots, t_k]}K$ is a compact subset of $R^k$. Trivially $\pi_{[t_1, t_2, \ldots, t_k]}^{-1} C_k\supset K$. Hence $\mu_n \pi_{[t_1, t_2, \ldots, t_k]}^{-1} C_k\geq \mu_n(K) > 1-\varepsilon$. This shows $(\mu_n \pi_{[t_1, t_2, \ldots, t_k]}^{-1})$ is a tight sequence of measures on $R^k$ and admits a weakly convergent subsequence.\\
Let $(\mu_{1, n}\;\pi_{t_1}^{-1}) $ be then a weakly convergent subsequence of $(\mu_n\;\pi_{[t_1]}^{-1}) $.  The tightness of the sequence $(\mu_{1, n})$
implies the tightness of the sequence $(\mu_{1, n}\;\pi_{[t_1, t_2]}^{-1})$ in $R^2$
which would then admit of a weakly convergent subsequence, say, $(\mu_{2, n}\;\pi_{[t_1, t_2]}^{-1})$. Now we start with the sequence $ (\mu_{2, n})$ and arguing similarly arrive at $(\mu_{3, n})$ which is such that $(\mu_{3, n}\;\pi_{[t_1]}^{-1}) $, $(\mu_{3, n}\;\pi_{[t_1, t_2]}^{-1})$ and $(\mu_{3, n}\;\pi_{[t_1, t_2, t_3]}^{-1})$ are weakly convergent sequences in $R^1, R^2, R^3$ respectively. In this way we
determine a family of sequences  $(\mu_{j, n}, n=1, 2, \ldots), j=1, 2, \ldots$ where $(\mu_{j+1, n})$ is a subsequence of $(\mu_{j, n})$. The
diagonal sequence $(\nu_n=\mu_{n, n})$ will have the property that $\nu_n \pi_{[t_1, t_2, \ldots, t_k]}^{-1}$  is a weakly convergent sequence for every
$k\geq 1$. It now follows by Theorem 2.10 that there exists a probability measure $\nu$ on $\mathcal{B}^*$ such that $\nu_n\;\pi_{[t_1, t_2, \ldots, t_k]}^{-1}$ converges weakly to  $\nu\;\pi_{[t_1, t_2, \ldots, t_k]}^{-1}$  for every $k \geq 1$. Because of this result and because of the fact that $(\nu_n)$ being a subsequence of  $(\mu_n)$ is
a tight sequence we conclude, by Theorem 2.9, that $\nu_n$ converges weakly to $\nu$.\\
\textbf{Remark.} It is of interest to know there are well defined and interesting metrics on $D$ which are regular. We now consider 3 such metrics on $D$.\\
(a)Let $\xi(x,y)=\sum_{k=1}^\infty\frac{1}{2^k}|x(t_k)-y(t_k)|+\sum_{k=1}^\infty\frac{1}{2^k}|x(t_k-)-y(t_k-)|$ where $0 \leq t_k\leq 1, k=1, 2, \ldots$. Note that this is well defined since members of $D$ are bounded functions.\\
(b) Let  $\rho$ be the uniform metric  on $D$.\\
(c) Let $d$ be the Skorohod metric on $D$.\\
It  can be shown (See Pakshirajan [6]) that these 3  metrics are regular.
\section{Prohorov space}
In this Section we define a Prohorov space and observe that a complete separable metric space is a Prohorov space. We discuss some examples of non-Prohorov spaces.\\
\textbf{Definition 3.1.}\\
A metric space $M$ is called a Prohorov space if every compact $\mathcal{F} \subset \mathcal{M}$, where $\mathcal{M}$ is the set of all  probability measures on $M$, is tight. \\
Preiss [9]  made the following two assertions concerning metric spaces in which Prohorov's theorem is not valid:\\
(a) A metric space $X$ which is of first category (see definition 3.3 below) in itself is not a Prohorov space, and\\
(b) If $X$ is a countable dense-in-itself metric space (e.g. the space of rational numbers) is not a Prohorov space.\\
\textbf{Definition 3.2.}\\
A subset of a topological space $X$ is said to be nowhere dense if its closure has  empty interior.\\
\textbf{Definition 3.3.}\\
A topological space $X$ is of first category in itself if it can be written as union of a countable number of sets $A_n \subset X$ such that each $A_n$ is  nowhere dense in $X$.\\
Now the set of rationals $\mathcal{Q}$, being a countable set, can be written as $\cup_{q \in \mathcal{Q}} \{q\}$ and is of first category in itself because $\{q\}$ is nowhere dense (also as a subset of $R$). Hence by (a) the set of rationals is not a Prohorov space.\\
If the metric space $M$ is complete and separable then by Theorem 2.5 it follows that $M$ is a Prohorov space. In this context we consider the space $M=\left\lbrace 0, 1, 2, \ldots\right\rbrace $ of non-negative integers. This is  a metric space with the distance function as the metric. This is complete and separable and hence  is a Prohorov space. It is of interest to note that  (a) is not applicable to $M$. It is not of first category in itself because $M= \cup_{n \in M}\; \left\lbrace n \right\rbrace $ where  $\left\lbrace n \right\rbrace$ is an open set in the topology induced by the usual metric. However as a subset of the real line, with usual topology, $M$ is of first category. \\
\noindent It has been of interest to find non-Prohorov spaces. Investigation of existence of a non-Prohorov space was initiated by Varadarajan [13] and his work enthused several people to work on non-Prohorov spaces. We refer to Preiss [9] and  Tops$\o$e [12] for other references and related problems. We now present  an example of a separable metric space which is \textbf{not} a Prohorov space, discussed by Tops$\o$e [12], with greater clarity. \\
\noindent Let $S = [0,\ 1] \times [0,\ 1] = I \times J$, say. Let $\pi$\ denote the projection operator from $S$\ on to $I$\ : If $(x,\ y) \in S$\ then $\pi(x,\ y) = x$. Define, for each $x \in I$, set $\Lambda(x) = \{y\ :\ y \in J,\ (x,\ y) \in S$, the section of $S$ at $x$. Let $\mathscr{K}$\ be the collection of all compact subsets $K$\ of $S$\ possessing the property $\pi(K) = [0, 1]$.\\
We shall first prove that there exists a  set $ A \subset S$, its projection under $\pi$ on the $ x-$ axis is a single point set  and it has nonempty intersection with every  compact set $K$ whose projection is the entire interval [0, 1].\\
\noindent Since a subset in $S$\ is compact iff it is closed, the cardinality of the family of compact sets is the same as that of the family of closed subsets, which is equal to the cardinality of the family of open sets. The cardinality of this last family is known to be $\mathfrak{c}$ , the cardinality of the continuum. Hence the cardinality of the collection of all compact sets is $\mathfrak{c}$. \\
\noindent Let $\mathscr{K}_1$\ consist of all compact sets $K$\ which are straight lines with end points on the lines  $x = 0\ \text{and}\ x = 1$. For example the set $\{(x,\ x)\ :\ x \in I\} \in \mathscr{K}_1$. The cardinality of $\mathscr{K}_1$\ is clearly $\mathfrak{c}$. \\
\noindent Hence $card{\mathscr{K}} = \mathfrak{c}$.\\
\noindent Since $\mathscr{K}$\ and $I$\ have the same cardinality (namely,\ $\mathfrak{c}$),  there exists a bijection $\Lambda$\ between $I$\ and $\mathscr{K}$\ (Schr$\ddot{\text{o}}$der-Bernstein theorem, see, p. 17, Kolmogorov and Fomin [5]). Given $x \in I,\ \Lambda(x)$\ will be the corresponding member in $\mathscr{K}$. Given $x$\ find $y = y(x)$\ such that $(x,\ y) \in \Lambda(x)$. This is possible since $\pi(\Lambda(x)) = I$.\\
\noindent Consider the set $A = \{(x,\ y(x)),\ x \in I\}$.  We note that if $K \in \mathscr{K}$, then there exists $u \in I$\ such that $\Lambda(u) = K$. This implies there exists $v$\ such that $(u,\ v) \in \textbf{K}$. Hence $(u,\ v) \in A$\ and we conclude\\ $A \cap K \neq \emptyset$.\hfill(3.1)\\
\noindent Define $\textbf{M} = A^{\prime}$. \\
We now claim that for every compact subset $\textbf{K}$  there exists an $x \in I$  such that the section of $K$ at $x$ is empty. \\ Note that $\textbf{M} = \bigcup\limits_{x \in I} \big( \pi^{-1}(\{x\}) \sim (x, y(x))\big)$. If $K$\ is a compact subset of \textbf{M}, then trivially $K \cap A = \emptyset$. Hence, by (3.1), $\pi(K) \neq [0,\ 1]$. This means that to every compact set $K \subset M$\ there can be found at least one $u \in I$\ such that\\ the line $x = u$\ has null intersection with \textbf{K}.\hfill(3.2)\\ 
We next find a family of probability measures on $\pi^{-1}(\{x\}), \; x\in I$ which is not tight.\\
\noindent Consider the space \textbf{M}\ endowed with the metric $d$\ inherited from $S$. Being a subset of a separable metric space, $(\textbf{M},\ d)$\ is separable (ref. p. 40,\ Theorem 1,  Zaanen [14] Lebesgue measure on $\pi^{-1}(\{x\}),\ x \in I$. This family considered extended to all of $S$\ is tight since $S$\ is a compact set. Since $S$\ is a complete and separable metric space it follows that it is sequentially compact (ref. Theorem 2.5).\\
Finally we find a family of probability measures on $\pi^{-1}(\{c\})\cap \textbf{M}$ which is a sequentially compact subset of the space of probability measures on $\textbf{M}$.\\
\noindent Let $\mu_x$\ denote the Lebesgue measure on $\pi^{-1}(\{x\}) \cap \textbf{M}$. i.e. on $Q_x = \pi^{-1}(x) \sim (x, y(x))$. Let $C$\ be an arbitrary closed subset of $Q_x$. There exists then a closed subset $C^*$\ of $\pi^{-1}(\{x\})$\ such that $C = C^* \cap \textbf{M}$. $C$\ will be equal to either $C^*$\ or $C^* \sim (x, y(x))$. As such $\mu_x(C) = \nu_x(C^*)$. This implies that  $(\nu_{x_n})$\ is weakly convergent iff $(\mu_{x_n})$\ is. We conclude that this family, $\mathcal{F}$,  of measures $\mu_x$\ on the separable  metric space $\pi^{-1}(\{x\}) \cap \textbf{M}$ is a sequentially compact subset of the space of probability measures on $M$. But it is not tight since for every compact subset $K \subset M$\ there exists a $\mu_x$\ with $\mu_x(K) = 0$\ (ref. (3.2))\\
With this the proof is complete that ($\textbf{M},\ d$)\ is a non-Prohorov space.\qed\\ 
\noindent \textbf{Remark.}\\
Space $(\textbf{M},\ d)$\ cannot be  a complete metric space since a complete and separable metric space is automatically a Prohorov space (ref. Theorem 2.5).

	\end{document}